\documentclass[12pt]{article}

 \usepackage{amsmath}
 \usepackage{amssymb}
\usepackage[margin=1.2 in]{geometry}
\usepackage{todonotes}
\usepackage{color}

 \renewcommand\footnotemark{}

\newtheorem{Theorem}{Theorem}[part]

\newtheorem{Proposition}[Theorem]{Proposition}

\newtheorem{Corollary}[Theorem]{Corollary}

\newtheorem{Remarks}[Theorem]{Remarks}

\makeatletter \@addtoreset{equation}{section}

\@addtoreset{Theorem}{section}

 \date{}
\font\tenmath=msbm10
\font\sevenmath=msbm7
\font\fivemath=msbm5
\newfam\mathfam \textfont\mathfam=\tenmath
\scriptfont\mathfam=\sevenmath
\scriptscriptfont\mathfam=\fivemath

\def \b{\noindent}
\def \e{I\!\!E}
\def \p{I\!\!P}

\def \={{\buildrel {\rm (law)} \over =}}
\def \N{{I\!\!N}}
\def \R{{I\!\!R}}

\def \1{1\!\!1}
\def \l{\ell}

\begin{document}
\title{L\'evy measures of infinitely divisible positive processes - examples and distributional identities}  
\author{Nathalie Eisenbaum \thanks{\hspace{-20pt} MAP5, CNRS - Universit\'e  de Paris;  nathalie.eisenbaum@parisdescartes.fr \hspace\fill} {and Jan Rosi\'nski\thanks{\hspace{-20pt} Department of Mathematics,  University of Tennessee; rosinski@math.utk.edu. \hspace\fill  \newline
 Rosinski's research was partially  supported by the Simons Foundation grant 281440.}}}

\date{}

\maketitle

\begin{abstract}
\noindent The law of a positive infinitely divisible process with no drift is  characterized by its L\'evy measure on the paths space.  Based on  recent results of the two authors, it is  shown that even for simple  examples of such process, the knowledge of their L\'evy measures  allows to obtain remarkable distributional identities.	
\end{abstract}

\section{Introduction}\label{s:intro} 

A random process  is infinitely divisible if all its finite dimensional marginals are infinitely divisible. 
Let $\psi = (\psi(x), x \in E)$ be a nonnegative infinitely divisible process with no drift.  The infinite divisibility of $\psi$ is characterized by the existence of a unique measure $\nu$ on $\R^E_+$, the space of all functions from $E$ into $\R_+$, such that for every $n>0$, every $\alpha_1, .., \alpha_n$ in $\R_+$ and every $x_1,..,x_n$ in $E$:
\begin{equation}\label{LaplaceT}
\e[ \exp\{ - \sum_{i = 1}^n \alpha_i \psi(x_i) \}] = 
\exp\{ - \int_{\R_+^E} (1 - e^{-\sum_{i = 1}^n \alpha_i y(x_i)}) \nu(dy) \}.
\end{equation}
The measure $\nu$ is called the L\'evy measure of $\psi$. 
The existence and uniqueness of such measures was established in complete generality in \cite{R}. In section \ref{prelim},  we recall some definitions and facts about L\'evy measures.  
%In particular,  if $\psi$ is as in \eqref{LaplaceT}, then $\int_{\R_+^E} f(x) \, \nu(df) < \infty$ for every $x \in E$, and $\psi$ admits the following representation:   
%\[
%   (\psi(x), x \in E) \, \=\, \Big( \int_{\R_+^E} f(x) \, N(df), \ x \in E\Big),
%\]
%where $N$ is a Poisson random measure with intensity $\nu$. 

It might be difficult to obtain an  expression for the L\'evy measure $\nu$ directly from (\ref{LaplaceT}).  In \cite{E}, a general expression for $\nu$ has been established. Its proof is based on several identities involving $\psi$.  Among them:

\noindent{\sl For every $a \in E$ with $ 0 < \e[\psi(a)] < \infty$, there exists a nonnegative process $(r^{(a)}(x) , x \in E)$ independent of $\psi$ such that}
\begin{equation}\label{IsoNat}
\psi + r^{(a)} \; \hbox{\sl has the law of} \; \psi \;  \hbox{\sl under} \; \e\Big[ {\psi(a) \over \e [\psi(a)]}, \; \cdot\; \Big]
\end{equation}
Actually, the existence of $(r^{(a)}, a \in E)$ characterizes the infinite divisibility of $\psi$. This characterization has been established in \cite{E1}, see also  \cite[Proposition 4.7]{R}.

Under an assumption of stochastic continuity for $\psi$, the general expression for $\nu$ obtained in \cite{E},  is the following: 
\begin{equation}\label{levyMeasure}
\nu(F) = \int_E \e\Big[ \frac{F(r^{(a)})}{\int_Er^{(a)}(x) m(dx)}\Big] \e[ \psi(a)] m(da) , 
\end{equation}
for any measurable functional $F$ on $\R_+^E$,
where $m$ is any $\sigma$-finite measure with support equal to $E$ such that $\int_E \e[ \psi(x)] m(dx) < \infty$.

Moreover the law of $r^{(a)}$ is connected to $\nu$ as follows (see \cite{E}, \cite{R}):
\begin{equation}\label{r}
\e[F( r^{(a)})] = \frac{1}{\e[\psi(a)]} \int_{\R_+^E} y(a) F(y)\, \nu(dy). 
\end{equation}

The problem of determining $\nu$ is hence equivalent to the one of the law of $r^{(a)}$ for every $a$ in $E$. But knowing $\nu$, one can not only write (\ref{IsoNat}) but many other identities of the same type. In each one, the process $r^{(a)}$ is replaced by a process  with an absolutely continuous law with respect  to $\nu$ (see \cite[Theorem 4.3(a)]{R}). 

Some conditionings on $\psi$ lead to a splitting of $\nu$. This allows to obtain  decompositions of $\psi$ into  independent infinitely divisible components (see \cite[Theorems 1.1, 1.2 and 1.3]{E}).  As an example:

\noindent{\sl For every $a \in E$,  there exists a nonnegative infinitely divisible process $({\cal L}^{(a)}(x), x \in E)$  independent of an infinitely divisible process $((\psi(x), x \in E) | \psi(a) = 0)$ such that}
\begin{equation}\label{Condition}
 \psi \=\ (\psi \ |\  \psi(a) = 0) \ + \ {\cal L}^{(a)}\, . 
\end{equation}

By \cite[Theorem 1.2]{E}, the processes $(\psi \,|\, \psi(a)=0)$ and ${\cal L}^{(a)}$ have the respective L\'evy measures $\nu_a$ and $\tilde{\nu}_a$, where
\begin{equation} \label{E1}
  \nu_{a}(dy)= \1_{\{y(a)=0 \}} \nu(dy)\ \textrm{and} \quad \tilde{\nu}_{a}(dy)= \1_{\{y(a)>0 \}} \nu(dy)\, .  
\end{equation}

\bigskip
 
 In section \ref{Illustrations}, to illustrate the relations and identities \eqref{LaplaceT}--\eqref{Condition} we choose to consider simple examples of nonnegative infinitely divisible processes. In each case the L\'evy measure is directly computable from (\ref{LaplaceT}) or from the stochastic integral representation of $\psi$ (see \cite{RR}). Thanks to (\ref{IsoNat}) and its extensions, and (\ref{Condition}), we present remarkable identities satisfied by the considered nonnegative infinitely divisible processes.  Moreover the general expression (\ref{levyMeasure})  provides alternative formulas for the L\'evy measure, which are also remarkable.  We treat the cases of Poisson  processes, Sato processes, stochastic convolutions and  tempered stable subordinators. We also point out a connection with  infinitely divisible random measures.  We end  section \ref{Illustrations} by reminding the case of infinitely divisible permanental processes which is the first case for which identities in law of the same type as  (\ref{IsoNat}) have been established.  In 
 this case,  such identities in law are called "isomorphism theorems" in reference  to the very first one established by Dynkin \cite{D} the so-called "Dynkin isomorphism Theorem".

When $\psi$ is an infinitely divisible permanental process, the two processes $r^{(a)}$ and ${\cal L}^{(a)}$ have the same law.   If moreover $\psi$ is a  squared Gaussian process, Marcus and Rosen \cite{MR} have established correspondences between path properties of $\psi$ and the ones of ${\cal L}^{(a)}$.  The  extension of these correspondences to general  infinitely divisible permanental processes has been undertaken by several authors (see \cite{EK}, \cite{E},  \cite{MR1} or \cite{MR2}). Similarly,  in section \ref{Correspondency},  we consider a general infinitely divisible nonnegative process $\psi$ and  state some  trajectories correspondences between $\psi$ and ${\cal L}^{(a)}$ resulting from an iteration of (\ref{Condition}) (see also \cite{R}). 
 
Finally, observing that given an infinitely divisible positive process $\psi$,  $r^{(a)}$ is not a priori  ``naturally" connected to $\psi$, we present, 
in section \ref{limit}, $r^{(a)}$ as the limit of a sequence of processes naturally connected to $\psi$.

\bigskip

\section{Preliminaries on L\'evy measures}\label{prelim}
In this section we recall some definitions and facts about general L\'evy measures given in \cite[Section 2]{R}. Some additional material can be found in \cite{R3}. Let $(\xi(x), x \in E)$ be a real-valued infinitely divisible process, where $E$ is an arbitrary  nonempty set. A measure $\nu$ defined on the cylindrical $\sigma$-algebra $\mathcal{R}^E$ of $\R^E$ is called the L\'evy measure of $\xi$ if the following two conditions hold: 
\begin{itemize}
  \item[(i)]  for every $x_1,\dots, x_n \in E$, the L\'evy measure of the random vector $(\xi(x_1), \dots, \xi(x_n))$ coincides with the  projection of $\nu$ onto $\R^{\{x_1,\dots,x_n\}}$, modulo the  mass at the origin;
  \item[(ii)]  $\nu(A) =\nu_{\ast}(A \setminus 0_E)$ for all  $A \in \mathcal{R}^E$, where $\nu_{\ast}$ denotes the inner measure and  $0_E$ is the origin of $\R^E$.
\end{itemize}
The L\'evy measure of an infinitely divisible process always exists and  (ii) guarantees its uniqueness. Condition (i) implies that $\int_{\R^E} (f(x)^2\wedge 1) \, \nu(df) < \infty$ for every  $x \in E$. 

A  L\'evy measure  $\nu$ is $\sigma$-finite if and only if  then there exists a countable set $E_0 \subset E$ such that 
\begin{equation}\label{sf}
  \nu\{f \in \R^{E}: f_{| E_0} = 0 \} =0.
\end{equation}
Actually, if (i) and \eqref{sf} hold, then does so (ii) and $\nu$ is a $\sigma$-finite L\'evy measure.

Condition \eqref{sf} is usually easy to verify. For instance, if an infinitely divisible process $(\xi(x), x \in E)$ is separable in probability, then its  L\'evy measure satisfies \eqref{sf}, so is $\sigma$-finite. The separability in probability is a weak assumption. It says that there is a countable set $E_1 \subset E$ such that for every $x \in E$ there is a sequence $(x_n) \subset E_1$ such that $\xi(x_n) \to \xi(x)$ in $\p$.  Infinitely divisible processes whose L\'evy measures do not satisfy \eqref{sf} include such pathological cases as an uncountable family of independent Poisson random variables with mean 1.

%Thus, $\nu$ is a $\sigma$-finite L\'evy measure of $\xi$ if and only if it satisfies (i) and \eqref{sf}. One may ask a question, how to ensure that $\nu$ is $\sigma$-finite if the process $\xi$ is known but its L\'evy measure is unknown? The necessary and sufficient condition for the $\sigma$-finiteness of $\nu$ in terms of $\xi$ is given in \cite[Theorem 2.14]{R} but the condition is quite involved. Fortunately, there is a simple and general sufficient condition which says that if an ID process $\xi$ is separable in probability then it has a $\sigma$-finite L\'evy measure. Recall that a process $(\xi(x), x \in E)$ is separable in probability if there is a countable $E_1 \subset E$ such that for every $x \in E$ there is a sequence $(x_n) \subset E_1$ such that $\xi(x_n) \to \xi(x)$ in $\p$.  

If the process $\xi$ has paths in some ``nice'' subspace of $\R^E$, then due to the transfer of regularity  \cite[Theorem 3.4]{R}, its L\'evy measure $\nu$ is carried by the same subspace of $\R^E$. Thus, one can investigate the canonical process on $(\R^E, \mathcal{R}^E)$  under the law of $\xi$ and also under the measure $\nu$, and relate their properties. This approach was successful in the study of distributional properties of subadditive functionals of paths of infinitely divisible processes \cite{RS} and the decomposition and classification of stationary stable processes \cite{R1}, among others.

If an infinitely divisible process $\xi$ without Gaussian component has the L\'evy measure $\nu$, then it can be represented as
\begin{equation}\label{L-I}
  (\xi(x), x \in E)\, \=\, \Big( \int_{\R^E} f(x) \, [N(df) - \chi(f(t)) \nu(df)] + b(x), \ x \in E\Big)
\end{equation}
where $N$ is a Poisson random measure on $(\R^E, \mathcal{R}^E)$ with intensity measure $\nu$, $\chi(u)= \1_{[-1,1]}(u)$, and $b \in \R^E$ is deterministic.  
Relation \eqref{L-I} can be strengthen to the equality almost surely under some minimal regularity conditions on the process $\xi$, provided the probability space is rich enough  (see \cite[Theorem 3.2]{R}). This is an extension to general infinitely divisible processes of the celebrated L\'evy-It\^o representation.

%If $\xi=(\xi_t, t \in \R_{+})$ is a L\'evy process without Gaussian part, then both $\xi$ and its L\'evy measure $\nu$ are carried by $D[0, \infty)$ and  \eqref{L-I} holds a.s. After a simple change of variable in the integral \eqref{L-I}, that formula is just the classical L\'evy-It\^o representation,  see Corollary 3.3 in \cite{R}.

Obviously, all the above apply to processes presented in the introduction but in more transparent form. Namely, if $(\psi(x), x \in E)$ is an infinitely divisible nonnegative process, then its L\'evy measure $\nu$ is concentrated on $\R_+^E$ and \eqref{L-I} becomes

\begin{equation}\label{NonnegativeCase}
 (\psi(x), x \in E) \, \=\, \Big(f_0(x) + \int_{\R_+^E} f(x) \, N(df), \ x \in E\Big),
\end{equation}
where $N$ is a Poisson random measure on $\R_+^E$ with intensity measure $\nu$ such that $\int_{\R_+^E} (f(x)\wedge 1) \, \nu(df) < \infty$ for every  $x \in E$.  Moreover, $\e [\psi(x)] < \infty$ if and only if  $\int_{\R_+^E} f(x) \, \nu(df) < \infty$ and $f_0\ge 0$ is a deterministic drift.

Since $N$ can be seen as a countable random subset of $\R_+^E$, one can write  (\ref{NonnegativeCase}) as 
\begin{equation}\label{NonnegativeCase2}
 (\psi(x), x \in E) \, \=\, \Big(f_0(x)+ \sum_{f \in N} f(x), \ x \in E\Big).
 \end{equation}

We end this section with a necessary and sufficient condition for a measure $\nu$ to be the L\'evy measure of a nonnegative infinitely divisible process. It is a direct consequence of \cite{R} section 2. From now on will assume that $\psi$ has no drift, in which case $f_0=0$ in \eqref{NonnegativeCase}-\eqref{NonnegativeCase2}.  

Let $\nu$ be a measure on
 $(\mathbb{R}_+^E, {\cal B}^E)$, where  ${\cal B}^E$ denotes the  cylindrical $\sigma$-algebra associated to $\mathbb{R}_+^E$ the space of all functions from $E$ into $\mathbb{R}_+$.  There exists an infinitely divisible nonnegative process $(\psi(x), x \in E)$ such that 
 for every $n > 0$, every $x_1,..,x_n$ in $E$:
\[
 \mathbb{E}[ \exp\{ - \sum_{i = 1}^n \alpha_i \psi(x_i) \}] = 
\exp\{ - \int_{\mathbb{R}_+^E} (1 - e^{-\sum_{i = 1}^n \alpha_i y(x_i)}) \nu(dy) \},
\]
if and only if $\nu$ satisfies the two following conditions:

(L1) for every $x \in E$   $\nu( y(x) \wedge 1) ) < \infty$, 

(L2) for every $A \in {\cal B}^E$, $\nu(A) = \nu_*(A \setminus 0_E)$, where $\nu_*$ is the inner measure.

\section{Illustrations}\label{Illustrations}

 By a standard uniform random variable we mean a random variable with the uniform law on $[0,1]$. A random variable with exponential law and mean 1 will  be called standard exponential. 

\subsection{Poisson process}

A Poisson process $(N_t, t \geq 0)$ with intensity $\lambda m$, where $\lambda > 0$ and $m$ is the Lebesgue measure on $\R_{+}$, is the simplest L\'evy process but its L\'evy measure $\nu$ is even simpler. It is a $\sigma$-finite measure given by
\begin{equation} \label{LM-P}
  \nu(F) = \lambda \int_0^{\infty} F\big(\1_{[s, \infty)}) \, ds, 
\end{equation}
for every measurable functional $F: \R_{+}^{[0,\infty)} \mapsto \R_{+}$.  Thus \eqref{LM-P} says that $\nu$ is the image of $\lambda m$ by the mapping $s \mapsto \1_{[s, \infty)}$ from $\R_{+}$ into $\R_{+}^{[0,\infty)}$. 

Formula \eqref{LM-P} is a special case of \cite[Example 2.23]{R}. We  will derive it here for the sake of illustration and completeness. 

Let  $(N_t, t \geq 0)$ be a Poisson process as above. 
By a routine computation of the Laplace transform, we obtain that  for every $0 \le t_1 < \dots < t_n$ the L\'evy measure $\nu_{t_1,\dots, t_n}$ of $(N_{t_1}, \dots, N_{t_n})$ is of the form
\[
\nu_{t_1,\dots, t_n} = \sum_{i=1}^{n} \lambda \Delta t_i  \, \delta_{\mathbf{u}_i} , 
\]
where $\Delta t_i=t_i-t_{i-1}$, $t_0=0$, and $\mathbf{u}_i=(\underbrace{0,\dots,0}_{ \ i-1 \, \text{times}}, 1,\dots, 1) \in \R^n$, $i=1,\dots,n$.

To verify that \eqref{LM-P} satisfies (i) of Section \ref{s:intro}, consider a finite dimensional functional  $F$, that is $F(f)= F_0(f(t_1), \dots, f(t_n))$, where $F_0: \R_{+}^n \mapsto \R_{+}$ is a Borel function with $F_0(0, \dots, 0)=0$ and $0 \le t_1 < \dots < t_n$. From \eqref{LM-P} we have
\begin{align*} 
   \nu(F) &= \lambda \int_0^{\infty} F\big(\1_{[s, \infty)}\big) \, ds = \lambda \int_0^{\infty} F_0(\1_{[s, \infty)}(t_1), \dots, \1_{[s, \infty)}(t_n)) \, ds \\
   &= \lambda \sum_{i=1}^{n} \int_{t_{i-1}}^{t_i} F_0(\mathbf{u}_i) \, ds = \int_{\R_{+}^{n}} F_0(x) \, \nu_{t_1,\dots, t_n}(dx)
\end{align*}
which proves (i). Condition \eqref{sf} holds for any unbounded set, for instance $E_0=\N$. Indeed,
\begin{align*} 
  \nu\{f\in \R_{+}^{[0, \infty)}: f_{| \N} =0\} = \lambda \int_0^{\infty} \1\{s: \1_{[s, \infty)}(n)=0 \ \forall n \in \N\})\, ds = 0,
\end{align*}
so that $\nu$ is the L\'evy measure of $(N_t, t \geq 0)$.

The next proposition exemplifies remarkable identities resulting from (\ref{Condition}) and \eqref{IsoNat}. It also gives an alternative ``probabilistic'' form of the L\'evy measure $\nu$.

\begin{Proposition}\label{Poisson} Let $N = (N_t, t \geq 0)$ be a Poisson process with intensity $\lambda m$, where $m$ is the Lebesgue measure on $\R_{+}$ and $\lambda>0$. 
\begin{itemize}
  \item[\rm{(a1)}] Given $a> 0$, let $r^{(a)}$ be the process defined by: $r^{(a)}(t) := \1_{[aU, \infty)}(t), \  t \ge 0$,  
  where $U$ is a standard uniform random variable independent of $(N_t, t \geq 0)$. Then $(r^{a}(t), t \geq 0)$  satisfies \eqref{IsoNat}, that is,
 \[
 (N_t \ + \  \1_{[aU,\, \infty)}(t) ,  \ t \geq 0)\, \=\, (N_t, \ t \geq 0) \ \hbox{under} \ \e\Big[\frac{N_a}{\lambda a}; \ .\ \Big]\, .
\]

  \item[\rm{(b1)}]For any nonnegative random variable $Y$ whose support equals $\R_+$ and $\e Y < \infty$,  the L\'evy measure $\nu$ of $(N_t, t \ge 0)$ can be represented as  
\[
 \nu(F) = \lambda \ \e\big[Y  h(UY) \ F\big(\1_{[UY,\, \infty)}\big) \big] 
  \]
   for every measurable functional $F: \R_{+}^{[0,\infty)} \mapsto \R_{+}$, 
where $U$ is a standard uniform random variable independent of $Y$ and  $h(x) =1/ \p [Y \geq x]$. 

In particular, if $Y$ is a standard exponential random variable independent of $U$, then
\[
 \nu(F) = \lambda \ \e\big[Y e^{UY} \ F\big(\1_{[UY,\, \infty)}\big) \big].
  \]

  \item[\rm{(c1)}]  The components of the decomposition \eqref{Condition}:  $N \= (N \,|\, N_a =0) + {\cal L}^{(a)}$,  \  can be identified as
\[
(N_t, t \ge 0 \,|\, N_a =0)\ \=\ (N_{t \vee a} - N_a, \ t \geq 0).
\]
and 
\[
({\cal L}^{(a)}_t ,  \ t \geq 0)\ \=\ (N_{t\wedge a}, \ t \geq 0).\]

The L\'evy measures $\nu_{a}$ and $\tilde{\nu}_{a}$ of  $(N_t, t \ge 0 \,|\, N_a=0)$  and of $({\cal L}^{(a)}_t, t \ge 0)$, respectively, are  given by
\[\nu_{a}(F) = \lambda \int_a^{\infty} F\big(\1_{[s, \infty)}\big) \, ds, \]
and
\[\tilde{\nu}_{a}(F) = \lambda \int_0^{a} F\big(\1_{[s, \infty)}\big) \, ds, \] 
  for every measurable functional $F: \R_{+}^{[0,\infty)} \mapsto \R_{+}$.

\end{itemize}

\end{Proposition}
\bigskip

\b {\bf Proof} \ 
(a1): By (\ref{r}) we have for any measurable functional $F: \R^{[0,\infty)} \mapsto \R_{+}$
\begin{align*} 
  \e F(r^{(a)}_t, t\ge 0) &= \frac{1}{\e N_a}\int F(y)\, y(a) \, \nu(dy)  \\
  &= \frac{1}{a} \int_0^{\infty} F(\1_{[s, \infty)})\, \1_{[s, \infty)}(a) \, ds \\
  &= \frac{1}{a} \int_0^{a} F(\1_{[s, \infty)}) \, ds = \e F(\1_{[aU, \infty)}).
\end{align*}
Thus $(r^{(a)}_t, t\ge 0) \=\, (\1_{[aU, \infty)}(t), t\ge 0)$. Choosing  $U$ independent of $N$, we have \eqref{IsoNat} for $r^{(a)}_t = \1_{[aU, \infty)}(t)$, $t\ge 0$, which completes the proof of (a1).

\smallskip

\b (b1): This point is an illustration of the invariance property in $m$ of (\ref{levyMeasure}).    Indeed, since the   process  $(N_t, t \geq 0)$ is  stochastically continuous   we have for every $\sigma$-finite measure $\tilde{m}$ whose support is $[0, \infty)$ and $\int_0^{\infty} t \, \tilde{m}(dt) < \infty$
\begin{align*} 
  \nu(F) &= \int_0^{\infty} \e\left[\frac{F(r^{(a)})}{\int_0^{\infty} r^{(a)}_{s} \, \tilde{m}(ds)}\right] \e[N_a]\, \tilde{m}(da) \smallskip \\
  &= \lambda \int_0^{\infty} \e\left[\frac{F(\1_{[aU, \infty)})}{\tilde{m}([aU, \infty))}\right] a \, \tilde{m}(da).
\end{align*}
If $\tilde{m}$ is the law of a nonnegative random variable $Y$, then
\begin{align*} 
   \nu(F) &= \lambda \int_0^{\infty} \e\left[a\, h(aU) F(\1_{[aU, \infty)}) \right]  \, \tilde{m}(da) \\ &= \lambda \e\left[Y h(UY) F(\1_{[UY, \infty)}) \right],
\end{align*}
which is the formula in (b1). 

%Implementing a formula established in  \cite{E} (Theorem 1.2) we get the L\'evy measure $\tilde{\nu}_a$ of ${\cal L}^{(a)}$: 
%$$ 
%\tilde{\nu}_a(F) = \e \left[\frac{\e[N_a]}{r^{(a)}_a} \1_{\{r^{(a)}_a>0 \}} F(r^{(a)}_t, \, t \ge 0)\right]. 
%$$
%Let $\tilde{F}(y)= \frac{1}{y(a)} \1_{\{y(a) >0\}}F(y)$, $y \in \R_{+}^{[0,\infty)}$. Then, from the first displayed formula in the proof of (b1), we have
%\begin{align*} 
%  \tilde{\nu}_a(F) &= \lambda a \, \e [\tilde{F}(r^{(a)}_t, \, t \ge 0)] = \int \1_{\{y(a)>0 \}}F(y(t), t\ge 0) \, \nu(dy)  \\
%&= \lambda \int_0^a F((\1_{[s, \infty)}(t), t\ge 0)) \, ds = \lambda a \e[ F(( \1_{[U_{a}, \infty)}(t), t \geq 0))] . 
%\end{align*}

(c1): Since $(N_t, t \ge 0 \,|\, N_a=0)$ has the L\'evy measure
$\nu_{a}(dy)= \1_{\{y(a)=0 \}} \nu(dy)$ (see \cite{E}), by \eqref{LM-P} we get
\begin{align*} 
  \nu_a(F) &= \int F(y) \, \1_{\{y(a)=0 \}}\nu(dy)  \\
&= \lambda \int_0^{\infty} F(\1_{[s, \infty)}) \1_{\{\1_{[s, \infty)}(a)=0 \}} \, ds\\
&= \lambda \int_a^{\infty} F(\1_{[s, \infty)})  \, ds\, .
\end{align*}
Since $\tilde{\nu}_a = \nu - \nu_a$, by \eqref{LM-P} we have
\[
\tilde{\nu}_a(F) = \lambda \int_0^{a} F(\1_{[s, \infty)}) 
 \, ds.
\]
 Let $0=t_0<t_1<\dots <t_n$ be such that  $t_m=a$ for some $m \le n$. For $\alpha_i>0$ we obtain
 \begin{align*} 
  \e \exp\Big\{-\sum_{i=1}^n \alpha_i & (\mathcal{L}^{(a)}_{t_i}-\mathcal{L}^{(a)}_{t_{i-1}}) \Big\} = \exp\{- \tilde{\nu}_a(1 - e^{-\sum_{i=1}^n \alpha_i (y(t_i)-y(t_{i-1}))}) \} \\
  &= \exp\{- \lambda \int_0^{a}(1 - e^{-\sum_{i=1}^n \alpha_i (\1_{[s,\infty)}(t_i)-\1_{[s,\infty)}(t_{i-1}))}) \,ds \} \\
  &= \exp\{- \lambda \sum_{i=1}^m\int_{t_{i-1}}^{t_i}(1 - e^{-\sum_{i=1}^n \alpha_i (\1_{[s,\infty)}(t_i)-\1_{[s,\infty)}(t_{i-1}))}) \,ds \} \\
  &= \exp\{- \lambda \sum_{i=1}^m (t_i - t_{t-1})(1 - e^{-\alpha_i})\} =   \e \exp\Big\{-\sum_{i=1}^n \alpha_i  (N_{t_i \wedge a}-N_{t_{i-1}\wedge a}) \Big\}
\end{align*}
which shows that $({\cal L}^{(a)}_t ,  \ t \geq 0)\, \=\, (N_{t\wedge a}, \ t \geq 0)$.  

Since $(N_{t\wedge a}, \ t \geq 0)$ and  $(N_{t \vee a} - N_a, \ t \geq 0)$    are independent and they add to $(N_t, t \ge 0)$,  $(N_t, t \ge 0 \,|\, N_a =0)\ \=\ (N_{t \vee a} - N_a, \ t \geq 0)$.  $\Box$

\begin{Remarks} {\ } \newline
{\rm
(1) By Proposition \ref{Poisson}(b1) the L\'evy measure $\nu$ of $N$ can be viewed as the law of the stochastic process
\[
(\1_{[UY, \infty)}(t),\, t \ge 0)
\]
under the infinite measure $\lambda Y h(UY)\, d\p$.  This point of view provides some intuition about the support of a L\'evy measure and better understanding how its  mass is distributed on the path space. 

(2) The process $(r^{(a)}_t, t\ge 0)$ of Proposition \ref{Poisson}(a1) is not infinitely divisible. Indeed, for each $t>0$, $r^{(a)}_t$ is a Bernoulli random variable.

(3) While the decomposition \eqref{Condition} is quite intuitive in case (c1), it is not so for general infinitely divisible random fields (cf. \cite{E}).

%
%(3) As mentioned in the introduction, we can generalize \eqref{IsoNat} into several potentially useful directions. For instance, from \cite[Example 4.11]{R} we have for any positive random variable $Z$ with density $g$ 
%\begin{equation} \label{}
%   (N_t \ + \  \1_{[Z, \infty)}(t) ,  \ t \geq 0)\, \=\, (N_t, \ t \geq 0) \ \hbox{under} \ \e\Big[\lambda^{-1} \int_0^{\infty} g(t) \, dN_t ; \ \cdot \ \Big].
%\end{equation}
%This gives Proposition \ref{Poisson}(b1) with $Z=U_a$ and $g(t)= a^{-1}\1_{[0, a]}(t)$.
%

%(4) In the present case, (\ref{Condition}) translates into the following identity, for every $a>0$
%\[ (N_t, t \geq 0) \ \= (N_{(t-a)^+}, t \geq 0) \ + \ (\bar{N}_{t \wedge a}, t \geq 0) \]
%where $\bar{N}$ is an independent copy of $N$. 

}	
\end{Remarks}
\bigskip

\subsection{Sato processes}

Recall that a process $X=(X_t,  t \ge 0 )$ is   $H$-self-similar, $H>0$, if for every $c>0$
$$
(X_{ct},  t \ge 0 ) \=\ (c^H X_t,   t \ge 0 )\, .
$$
It is well-known that a L\'evy process is $H$-self-similar if and only it is strictly $\alpha$-stable with $\alpha = 1/H \in (0,2]$, see \cite[Proposition 13.5]{S1}. In short, there are only obvious examples of self-similar L\'evy  processes.

Sato \cite{S} showed that within a larger class of additive processes there is a rich family of  self-similar processes which is generated by selfdecomposable laws. These processes are known as {\it Sato processes} and will be precisely defined below.

Recall that the law of a random variable $S$ is said to be selfdecomposable  if for every $b >1$ there exists an independent of $S$ random variable $R_b$ such that 
$$S \= b^{-1}S + R_b.$$
Wolfe \cite{W} and Jurek and Vervaat \cite{JV}, showed that a random variable $S$ is selfdecomposable if and only if 
\begin{equation} \label{JV}
  S \= \int_{0}^{\infty} e^{-s} \, dY_{s}
\end{equation}
for some L\'evy process $Y=(Y_s, \, s\ge 0)$ with $\e (\ln^{+} |Y_1|) < \infty$. Moreover, there is a 1-1 correspondence between the laws of $S$ and $Y_1$.  The process $Y$ is called the background driving L\'evy process (BDLP) of $S$. 

Sato \cite{S} proved that a random variable $S$ has the  selfdecomposable law if and only if  for each $H>0$ there exists a unique additive $H$-self-similar process $(X_t,  t \ge 0 )$ such that $X_1 \= S$. 
An additive self-similar processes, whose law at time 1 is selfdecomposable, will be called a Sato processes.

Jeanblanc, Pitman, and Yor \cite[Theorem 1]{JPY} gave the following representation of  Sato processes. Let $Y$ be the BDLP specified in \eqref{JV} and let $\hat{Y}=(\hat{Y}_s, \, s\ge 0)$ be an independent copy of $Y$. Then, for each $H>0$, the process
\begin{equation} \label{JPY}
  X_r := \begin{cases}
    \int_{\ln(r^{-1})}^{\infty}  e^{-Ht} \, d_t(Y_{Ht}) &  \text{if } 0\le r \le 1  \bigskip \\ 
   X_1  +  \int_{0}^{\ln r}  e^{Ht} \, d_t(\hat{Y}_{Ht})     & \text{if } r \ge 1. 
\end{cases}
\end{equation}
is the Sato process with selfsimilarity exponent $H$. Stochastic integrals in \eqref{JV} and \eqref{JPY} can be evaluated pathwise by parts due to the smoothness of the integrants. We will give another form of this representation that is easier to use for our purposes.

\begin{Theorem}\label{t:JR}
Let $\bar{Y}=(\bar{Y}_s, \, s\in \R)$ be a double sided L\'evy process such that $\bar{Y}_0=0$
and $\e (\ln^{+} |\bar{Y}_{1}|) < \infty$. Then, for each $H>0$, the process
\begin{equation} \label{JR1}
  X_t := \int_{\ln(t^{-H})}^{\infty} e^{-s} \, d\bar{Y}_s\, , \quad t \ge 0,
\end{equation}
is a Sato process with selfsimilarity exponent $H$. Conversely, any Sato process with selfsimilarity exponent $H$ has a version given by \eqref{JR1}.
\end{Theorem}

{\bf Proof } By definition, a double sided L\'evy process $\bar{Y}$ is indexed by $\R$, has stationary and independent increments, c{\`a}dl{\`a}g paths, and $\bar{Y}_0=0$ a.s. Since  \eqref{JR1} coincides with  \eqref{JV} when $t=1$, the improper integral $X_1=\int_{0}^{\infty} e^{-s} \, d\bar{Y}_s$ converges a.s. and it has a selfdecomposable distribution. Moreover, 
$$
X_{0+} = \lim_{t \downarrow 0} \int_{\ln(t^{-H})}^{\infty} e^{-s} \, d\bar{Y}_s = 0 \quad a.s.
$$

For every $0 < t_1 < \cdots < t_n$ and $u_k = \ln(t_k^{-H})$ the increments
$$
X_{t_k} - X_{t_{k-1}} = \int_{u_k}^{\infty} e^{-s} \, d\bar{Y}_s - \int_{u_{k-1}}^{\infty} e^{-s} \, d\bar{Y}_s = 
\int_{u_{k}}^{u_{k-1}} e^{-s} \, d\bar{Y}_s, \quad k=2,\dots, n
$$
are independent as $\bar{Y}$ has independent increments. Thus $X$ is an additive process. 

To prove the $H$-selfsimilarity of $X$, notice that since $X$ is an additive process, it is enough to show that for every $c>0$ and $0< t < u$ 
\begin{equation} \label{JR3}
  X_{cu} - X_{ct} \= c^H(X_{u} - X_{t}). 
\end{equation}
Since $\bar{Y}$ has stationary increments, we get
\begin{align*} 
  X_{cu} - X_{ct} &= \int_{\ln((cu)^{-H})}^{\ln((ct)^{-H})} e^{-s} \, d\bar{Y}_s = \int_{\ln(u^{-H}) + \ln(c^{-H})}^{\ln(t^{-H}) + \ln(c^{-H})} e^{-s} \, d\bar{Y}_s \bigskip\\
  & \= \int_{\ln(u^{-H})}^{\ln(t^{-H})} e^{-s - \ln(c^{-H})} \, d\bar{Y}_s = c^{H} (X_{u} - X_{t})\, ,
\end{align*}
which proves \eqref{JR3}. 

Conversely, let $X=(X_t: t \ge 0)$ be a $H$-selfsimilar Sato process. By \eqref{JV} there exists a unique in law L\'evy process $Y=(Y_t: t \ge 0)$ such that $\e (\ln^{+} |Y_1|) < \infty$ and
$$
X_1 \=\ \int_{0}^{\infty} e^{-s} \, dY_s\, .
$$
Let $Y^{(1)}$ and $Y^{(2)}$ be  independent copies of the L\'evy process $Y$. Define $\bar{Y}_s = Y^{(1)}_s$ for $s\ge 0$ and $\bar{Y}_s = Y^{(2)}_{(-s)-}$ \ for $s < 0$. Then $\bar{Y}$ is a double sided L\'evy process with $\bar{Y}_1 \=\ Y_1$.  
Then
$$
\tilde{X}_t := \int_{\ln(t^{-H})}^{\infty} e^{-s} \, d\bar{Y}_s\, , \quad t \ge 0,
$$
is a version of $X$.
$\Box$

\begin{Corollary}\label{c:JR2}
Let $X=(X_t: t \ge 0)$ be a $H$-selfsimilar Sato process given by \eqref{JR1}. Let $\rho$ be the L\'evy measure of $\bar{Y}_1$. 
Then the L\'evy measure $\nu$ of $X$ is given by
\begin{equation} \label{LM-S}
  \nu(F) = \int_{\R}\hspace{-2pt}\int_{\R} F(x e^{-s} \1_{[e^{-s/H}, \infty)}) \, \rho(dx) ds.
\end{equation}
\end{Corollary}
{\bf Proof} We can write \eqref{JR1} as $X_t= \int_{\R} f_t(s) \, d\bar{Y}_s$, where $f_t(s) = e^{-s} \1_{[e^{-s/H}, \infty)}(t)$. It follows from \cite[Theorem 2.7(iv)]{RR} that the L\'evy measure $\nu$ of $X$ is the image of $m \otimes \rho$ by the map $(s, x) \mapsto x f_{(\cdot)}(s)$ from $\R^{2}$ into $\R^{[0, \infty)}$. $\Box$
\bigskip

From now on we will consider a {\it $H$-selfsimilar nonnegative Sato process with finite mean and no drift}. By Theorem \ref{t:JR} we have
\begin{equation} \label{JRs1}
  \psi(t) = \int_{\ln(t^{-H})}^{\infty} e^{-s} \, d\bar{Y}_s\, , \quad t \ge 0,
\end{equation}
where $\bar{Y}=(\bar{Y}_t, t \in \R)$ is a double sided subordinator without drift such that $\bar{Y}_0=0$ and $\e \bar{Y}_{1} < \infty$. Consequently, $\e \psi(t) = \kappa t^{H}$, $t \ge 0$, where  $\kappa :=\e \psi(1) = \e \bar{Y}_1$.

\begin{Proposition}\label{p:JR2}

	Let $(\psi(t), t \ge 0)$ be a nonnegative $H$-selfsimilar Sato process given by \eqref{JRs1}. Therefore, the L\'evy measure $\rho$ of $\bar{Y}_1$ is concentrated on $\R_{+}$. 
	\begin{itemize}
  \item[\rm{(a2)}] Given $a> 0$, let $(r^{(a)}(t), t \geq 0)$ be the process defined by:
  \[r^{(a)}(t):= a^{H} UV \1_{[aU^{1/H}, \infty)}(t),  t \ge 0,\]
 where $U$ is a standard uniform random variable and  $V$ has the distribution $\kappa^{-1} x\rho(dx)$, with $U, V$ and $(\psi(t), t \ge 0)$  independent. Then $r^{(a)}$ satisfies \eqref{IsoNat}, that is, 
  $$
  \{\psi(t) + a^{H} UV \1_{[aU^{1/H}, \infty)}(t), \ t \ge 0 \} \=\ \{\psi(t), \ t \ge 0 \} \ \hbox{\rm{under}} \ \e\Big[\frac{\psi(a)} {\kappa a^{H}} ; \ .\ \Big]\, .
  $$
 \item[\rm{(b2)}]  Let $G$ be a standard exponential random variable, $U$ and $V$ be as above, and assume that $G$, $U$, and $V$ are independent.  Then the L\'evy measure $\nu$ of the process $(\psi(t), t \ge 0)$ can be represented as 
 \[
\nu(F) = \kappa \e\left[(UV)^{-1} e^{GU^{1/H}} F(G^{H}UV  \1_{[GU^{1/H}, \infty)}) \right]\,  \]
for every measurable functional $F: \R_{+}^{[0,\infty)} \mapsto \R_{+}$. Therefore, $\nu$ is the law of the process $(G^{H}UV  \1_{[GU^{1/H}, \infty)}(t), \,  t\ge 0)$ under the measure $ \kappa (UV)^{-1} e^{GU^{1/H}} \, d\p$.

\item[\rm{(c2)}] The components of the decomposition \eqref{Condition}:  $\psi \= (\psi \,|\, \psi(a) =0) + {\cal L}^{(a)}$, \  can be identified as
\[
(\psi(t), t \ge 0 \,|\, \psi(a) =0)\ \=\ (\psi(t \vee a) -\psi(a), \ t \geq 0).
\]
and 
\[
({\cal L}^{(a)}_t ,  \ t \geq 0)\ \=\ (\psi(t\wedge a), \ t \geq 0).\]

The L\'evy measures $\nu_{a}$ and $\tilde{\nu}_{a}$ of  $(\psi(t), t \ge 0 \,|\, \psi(a)=0)$  and of $({\cal L}^{(a)}_t, t \ge 0)$, respectively, are  given by
\[\nu_{a}(F) =  \int_{-\infty}^{\ln (a^{-H})}\hspace{-2pt}\int_{0}^{\infty} F(x e^{-s} \1_{[e^{-s/H}, \infty)}) \, \rho(dx) ds\, 
\]
and
\[\tilde{\nu}_{a}(F) = \int_{\ln (a^{-H})}^{\infty}\hspace{-2pt}\int_{0}^{\infty} F(x e^{-s} \1_{[e^{-s/H}, \infty)}) \, \rho(dx) ds, \] 
  for every measurable functional $F: \R_{+}^{[0,\infty)} \mapsto \R_{+}$. 
\end{itemize}
\end{Proposition}
 
{\bf Proof } \ (a2): By (\ref{r}) we have for any measurable functional $F: \R^{[0,\infty)} \mapsto \R_{+}$
\begin{align*} 
  \e F(r^{a}_t, \ t\ge 0) &= \frac{1}{\e \psi(a)} \int_{\R^{E}_{+}} F(y) y(a) \, \nu(dy) \\  & = \frac{1}{a^H\e \psi(1)}\int_{\R}\int_{\R_+} F(xe^{-s} \1_{[e^{-s/H}, \infty)})\, x e^{-s} \1_{[e^{-s/H}, \infty)}(a)\, \rho(dx) ds  \\
  & = \frac{a^{-H}}{\e \psi(1)} \int_{\ln (a^{-H})}^{\infty}\int_{\R_+} F(xe^{-s} \1_{[e^{-s/H}, \infty)})\,\, x  \, \rho(dx)\, e^{-s} ds \\
    & = a^{-H} \int_{\ln (a^{-H})}^{\infty} \e F(V e^{-s} \1_{[e^{-s/H}, \infty)})\, e^{-s} ds \\
  & = \e\Big[ F(a^HUV \1_{[aU^{1/H}, \infty)})\Big]\, .
\end{align*}
Thus $(r^{a}_t, \ t\ge 0) \= (a^HUV \1_{[aU^{1/H}, \infty)}(t), \ t\ge 0)$.
Since  $U$, $V$ and $\psi$ are independent, \eqref{IsoNat}  completes the proof of (a2).

(b2): Since the   process  $(\psi(t), t \geq 0)$ is  stochastically continuous we have for every $\sigma$-finite measure $\tilde{m}$ whose support is $[0, \infty)$ and $\int_0^{\infty} t^H \, \tilde{m}(dt) < \infty$
\begin{align*} 
  \nu(F) &= \int_0^{\infty} \e\left[\frac{F(r^{(a)})}{\int_0^{\infty} r^{(a)}_{s} \, \tilde{m}(ds)}\right] \e[\psi(a)]\, \tilde{m}(da) \smallskip \\
  &= \e [\psi(1)] \int_0^{\infty} \e\left[\frac{F(a^H UV\1_{[aU^{1/H}, \infty)})}{UV\, \tilde{m}([aU^{1/H}, \infty))}\right]  \, \tilde{m}(da).
\end{align*}
If $\tilde{m}$ is the law of a nonnegative random variable $W$, then
\begin{align*} 
   \nu(F) &= \e [\psi(1)] \int_0^{\infty} \e\left[\frac{h(aU^{1/H})}{UV}\,  F(a^H UV \1_{[aU^{1/H}, \infty)}) \right]  \, \tilde{m}(da) \\ 
   &=  \e [\psi(1)] \e\left[\frac{h(U^{1/H}W)}{UV} \,  F(UVW^H  \1_{[U^{1/H}W, \infty)}) \right]\,
\end{align*}
which is the formula in (b2). 

(c2): Since the conditional process $(\psi(t), t \ge 0 \,|\, \psi(a)=0)$ has the L\'evy measure $\nu_{a}(dy)= \1_{\{y(a)=0 \}} \nu(dy)$ (see \cite{E}), by \eqref{LM-S} we obtain for any  measurable functional $F: \R_{+}^{[0,\infty)} \mapsto \R_{+}$ and $a>0$
\begin{align*} 
  \nu_a(F) &= \int F(y) \, \1_{\{y(a)=0 \}}\nu(dy)  \\
  &= \int_{\R}\hspace{-2pt}\int_{\R_{+}} F(x e^{-s} \1_{[e^{-s/H}, \infty)})\1_{\{x e^{-s} \1_{[e^{-s/H}, \infty)}(a)=0 \}} \, \rho(dx) ds\\
&= \int_{-\infty}^{\ln (a^{-H})}\hspace{-2pt}\int_{0}^{\infty} F(x e^{-s} \1_{[e^{-s/H}, \infty)}) \, \rho(dx) ds\, .
\end{align*}

Since $\tilde{\nu}_a = \nu - \nu_a$,
$$
\tilde{\nu}_a(F) = \int_{\ln (a^{-H})}^{\infty}\hspace{-2pt}\int_{0}^{\infty} F(x e^{-s} \1_{[e^{-s/H}, \infty)}) \, \rho(dx) ds
$$
Let $0=t_0<t_1<\dots <t_n$ be such that  $t_m=a$ for some $m \le n$. For $\alpha_i>0$ we obtain
 \begin{align*} 
  \e & \exp\Big\{-\sum_{i=1}^n \alpha_i  (\mathcal{L}^{(a)}_{t_i}-\mathcal{L}^{(a)}_{t_{i-1}}) \Big\} = \exp\{- \tilde{\nu}_a(1 - e^{-\sum_{i=1}^n \alpha_i (y(t_i)-y(t_{i-1}))}) \} \\
  &= \exp\Big\{-  \int_{\ln (a^{-H})}^{\infty}\hspace{-2pt}\int_{0}^{\infty}(1 - e^{-\sum_{i=1}^n \alpha_i x e^{-s}(\1_{[e^{-s/H},\infty)}(t_i)-\1_{[e^{-s/H},\infty)}(t_{i-1}))}) \,\rho(dx) ds \Big\} \\
  &= \exp\Big\{-  \int_{\ln (a^{-H})}^{\infty}\hspace{-2pt}\int_{0}^{\infty}(1 - e^{-\sum_{i=1}^n \alpha_i x e^{-s}\1_{(t_{i-1}, t_i]}(e^{-s/H})}) \,\rho(dx) ds \Big\} \\
  &= \exp\Big\{- \sum_{i=1}^m\int_{\ln (t_i^{-H})}^{\ln (t_{i-1}^{-H})}\hspace{-2pt}\int_{0}^{\infty}(1 - e^{-\alpha_i x e^{-s}}) \,\rho(dx) ds \Big\} = 
  \prod_{i=1}^m \e \exp\big\{-\alpha_i(\psi(t_i) - \psi(t_{i-1}))\big\} \\
  &=\e \exp\{- \sum_{i=1}^n \alpha_i(\psi(t_i \wedge a)-\psi(t_{i-1}\wedge a))\big\}\, ,
\end{align*}
which shows that $({\cal L}^{(a)}_t ,  \ t \geq 0)\, \=\, (\psi(t\wedge a), \ t \geq 0)$. 

Since $(\psi(t \wedge a), \ t \geq 0)$ and $(\psi(t \vee a) - \psi(a), \ t \geq 0)$ are independent and they add to $(\psi(t),\, t \ge 0)$, we get $(\psi(t), t \ge 0 \,|\, \psi(a) =0)\ \=\ (\psi(t \vee a) - \psi(a), \ t \geq 0)$.  $\Box$

\subsection{Stochastic convolution}

Let $Z=(Z_t, t \ge 0)$ be a subordinator with no drift. For a fixed function $f: \R_{+} \mapsto \R_{+}$ and $t \ge 0$, the stochastic convolution  $f \ast Z$ is given  by
$$
 (f \ast Z)(t) = \int_{0}^{t} f(t-s) \, dZ_s\, .
 $$
Assume  that $\kappa := \e Z_1 \in (0, \infty)$ and $\int_{0}^t f(s) \, ds <\infty$ for every $t>0$. Therefore,  $\e [(f \ast Z)(t)] = \kappa \int_{0}^t f(s) \, ds < \infty$.  Set $f(u)=0$ when $u<0$.

We will consider the stochastic convolution process
\begin{equation}\label{sc}
  \psi(t) := \int_{0}^{t} f(t-s) \, dZ_s\, , \quad t \ge 0.
\end{equation}

Clearly, $(\psi(t),\ t \ge 0)$ is an infinitely divisible process. To determine its L\'evy measure  we write $\psi(t)= \int_{0}^{\infty} f_t(s) \, dZ_s$, where $f_t(s) = f(t-s)$.  It follows from \cite[Theorem 2.7(iv)]{RR} that the L\'evy measure $\nu$ of the process $\psi$ is the image of $m \otimes \rho$ by the map $(s, x) \mapsto x f_{(\cdot)}(s)$ acting from $\R_{+}^{2}$ into $\R^{[0, \infty)}$. That is,
\begin{equation} \label{LM-Conv}
   \nu(F) = \int_0^{\infty}\hspace{-6pt}\int_0^{\infty} F(x f(t-s), \, t \ge 0) \, \rho(dx) ds
\end{equation}

for every measurable functional $F: \R_{+}^{[0,\infty)} \mapsto \R_{+}$.

\begin{Proposition}\label{Convolution} Let $(\psi(t),\ t \ge 0)$ be a stochastic convolution  process as in \eqref{sc}. Let $\rho$ be the L\'evy measure of $Z_1$ and $I(a) :=  \int_{0}^a f(s) \, ds$.    
\begin{itemize}
  \item[\rm{(a3)}] Given $a> 0$ such that $I(a)>0$, let $r^{(a)}$ be the process defined by:
  \[r^{(a)}(t):= V f(t-U_a),  \quad t \geq 0\]
where the random variable $U_a$ has density $\dfrac{f(a-s)}{I(a)}$ on $[0,a]$, $V$ has the law $\kappa^{-1} x\rho(dx)$ on $\R_{+}$, and  $U_a$, $V$, and $(\psi(t): t \ge 0)$ are independent. 
Then $r^{(a)}$ satisfies \eqref{IsoNat}, that is,  
 \[
 \big(\psi(t) \ + \  V f(t-U_a),  \, t \geq 0\big)\, \=\, \big(\psi(t), \ t \geq 0\big) \ \hbox{under} \ \e\Big[\frac{\psi(a)}{\kappa I(a)}; \ .\ \Big]
\]

\item[\rm{(b3)}] Suppose that $\int_{0}^{\infty} e^{-\theta s} f(s)\, ds < \infty$ for some $\theta>0$.  Let $Y$ be a random variable with the exponential law of mean $\theta^{-1}$ and independent of $V$ specified in (a3).  Then the L\'evy measure $\nu$ of $(\psi(t), t \ge 0)$ can be represented as  
\[
 \nu(F) = \frac{\kappa}{\theta}  \e\left[V^{-1}e^{\theta Y} F(V f(t-Y), \,  t\ge 0)\right]. 
  \]
   for every measurable functional $F: \R_{+}^{[0,\infty)} \mapsto \R_{+}$. Therefore, $\nu$ is the law of the process $(V f(t-Y), \,  t\ge 0)$ under the measure $ \kappa \theta^{-1} V^{-1} e^{\theta Y}  \, d\p$.

\item[\rm{(c3)}] The components of the decomposition \eqref{Condition}:  $\psi \= (\psi \,|\, \psi(a) =0) + {\cal L}^{(a)}$, \  can be identified as
\[
(\psi(t), t \ge 0 \,|\, \psi(a) =0)\ \=\ \Big(\int_{0}^{t} f(t-s) \1_{D_a}(s) \, dZ_s, \ t \geq 0\Big)
\]
  and 
\[
({\cal L}^{(a)}_t ,  \ t \geq 0)\ \=\ \Big(\int_{0}^{t} f(t-s) \1_{D_a^{c}}(s) \, dZ_s, \ t \geq 0\Big)\]
where $D_a = \{s\ge 0: f(a-s)=0 \}$ and $D_a^{c} = \R_{+}\setminus D_a$\, .

The L\'evy measures $\nu_{a}$ and $\tilde{\nu}_{a}$ of  $(\psi(t), t \ge 0 \,|\, \psi(a)=0)$  and of $({\cal L}^{(a)}_t, t \ge 0)$, respectively, are  given by
\[\nu_{a}(F) =  \int_{D_a}\hspace{-2pt}\int_{0}^{\infty} F(x f(t-s), \, t \ge 0) \, \rho(dx) ds \, 
\]
and
\[\tilde{\nu}_{a}(F) = \int_{D_a^{c}}\hspace{-1.5pt}\int_{0}^{\infty} F(x f(t-s), \, t \ge 0) \, \rho(dx) ds, \] 
  for every measurable functional $F: \R_{+}^{[0,\infty)} \mapsto \R_{+}$. 
\end{itemize}
	
\end{Proposition}

{\bf Proof} 
(a3): From (\ref{r}) and \eqref{LM-Conv} we get 
\begin{align*} 
  \e F(r^{(a)}_t, \, t\ge 0) &= \frac{1}{\e \psi(a)}\int F(y)\, y(a) \, \nu(dy)  \\
  &= \frac{1}{\kappa I(a)} \int_0^{\infty}\hspace{-4pt}\int_0^{\infty} F\big(xf(t-s), \, t\ge 0\big)\, xf(a-s) \, \rho(dx) ds \\
  &=  \int_0^{a}\hspace{-4pt}\int_0^{\infty} F(x f(t-s), t\ge 0)\,   \, \frac{x \rho(dx)}{\kappa} \frac{f(a-s) ds}{I(a)} \\
  & =\e \big[F(V f(t-U_a), \,  t\ge 0)\big]\, .
\end{align*}

(b3): Since $\psi$ is stochastically continuous, using (\ref{levyMeasure}) and (a3), we have for every $\sigma$-finite measure $\tilde{m}$ whose support is $[0, \infty)$ and $\int_{0}^{\infty} I(a) \, \tilde{m}(da)< \infty$
\begin{align*} 
  \nu(F) &= \int_0^{\infty} \e\left[\frac{F(r^{(a)})}{\int_0^{\infty} r^{(a)}_{s} \, \tilde{m}(ds)}\right] \e[\psi(a)]\, \tilde{m}(da) \\
  &= \kappa \int_0^{\infty} \e\left[\frac{F(V f(t-U_a), \,  t\ge 0)}{V\int_0^{\infty} f(s-U_a) \, \tilde{m}(ds)}\right] I(a) \, \tilde{m}(da).
\end{align*}
Since $\tilde{m}$ is the law of $Y$ in our case, it is easy to check that $\beta:=\int_{0}^{\infty} I(a) \, \tilde{m}(da)< \infty$. Also,
$$
\int_0^{\infty} f(s-U_a) \, \tilde{m}(ds) = \beta \theta e^{-\theta U_a}\, .
$$
Then we get 
\begin{align*} 
   \nu(F) &= \frac{\kappa}{\beta \theta}\int_0^{\infty} \e\left[V^{-1} e^{\theta U_a} F(V f(t-U_a), \,  t\ge 0)\right] I(a) \theta e^{-\theta a} \, da \\
   &= \frac{\kappa}{\beta \theta}\int_0^{\infty} \int_0^a  \e\left[V^{-1} e^{\theta s}F(V f(t-s), \,  t\ge 0)\right] f(a-s) \, ds \, \theta e^{-\theta a} \, da \\
    &= \frac{\kappa}{\theta} \int_0^{\infty} \e\left[V^{-1} e^{\theta s} F(V f(t-s), \,  t\ge 0)\right]  \theta e^{-\theta s}\, ds \\
    &= \frac{\kappa}{\theta}  \e\left[V^{-1}e^{\theta Y} F(V f(t-Y), \,  t\ge 0)\right].
    \end{align*}

(c3): Since the conditional process $(\psi(t), t \ge 0 \,|\, \psi(a)=0)$ has the L\'evy measure $\nu_{a}(dy)= \1_{\{y(a)=0 \}} \nu(dy)$ (see \cite{E}), by \eqref{LM-S} we obtain for any  measurable functional $F: \R_{+}^{[0,\infty)} \mapsto \R_{+}$ and $a>0$
\begin{align*} 
  \nu_a(F) &= \int F(y) \, \1_{\{y(a)=0 \}}\nu(dy)  \\
  &= \int_0^{\infty}\hspace{-6pt}\int_0^{\infty} F(x f(t-s), \, t \ge 0) \, \1_{\{(x, s): \, xf(a -s))=0 \}} \, \rho(dx) ds\\
    &= \int_0^{\infty}\hspace{-6pt}\int_0^{\infty} F(x f(t-s), \, t \ge 0) \, \1_{D_a}(s) \, \rho(dx) ds\, .
\end{align*}
Using again \cite[Theorem 2.7(iv)]{RR} we see that $\nu_a$ is the L\'evy measure of the process 
$$\Big(\int_{0}^{t} f(t-s) \1_{D_a}(s) \, dZ_s, \ t \geq 0\Big)$$
which is a nonnegative infinitely divisible process without drift. Since the law of such process is completely characterized by its L\'evy measure, we infer that 
\[
(\psi(t), t \ge 0 \,|\, \psi(a) =0)\ \=\ \Big(\int_{0}^{t} f(t-s) \1_{D_a}(s) \, dZ_s, \ t \geq 0\Big)\, .
\]
Since $\tilde{\nu}_a = \nu - \nu_a$ and $\psi \= (\psi \,|\, \psi(a) =0) + {\cal L}^{(a)}$, we can apply the same argument as above to get 
\[({\cal L}^{(a)}_t ,  \ t \geq 0)\ \=\ \Big(\int_{0}^{t} f(t-s) \1_{D_a}(s) \, dZ_s, \ t \geq 0\Big)\, .
\]

$\Box$

%Implementing a formula from \cite{E} for the L\'evy measure $\tilde{\nu}_a$ of ${\cal L}^{(a)}$ we get
%$$ 
%\tilde{\nu}_a(H) = \e \left[\frac{\e[\psi(a)]}{r^{(a)}_a} \1_{\{r^{(a)}_a>0 \}} H(r^{(a)}_t, \, t \ge 0)\right]. 
%$$
%Let $\tilde{H}(y)= \frac{1}{y(a)} \1_{\{y(a) >0\}}H(y)$, $y \in \R_{+}^{[0,\infty)}$. Then, from (a2) we get
%\begin{align*} 
%  \tilde{\nu}_a(H) &= \kappa F(a) \, \e [\tilde{H}(r^{(a)}_t, \, t \ge 0)] \\
%  &= \kappa F(a) \, \e \Big[\frac{1}{W f(a-V_a)} \1_{\{V_a < a \}} H(W f(t-V_a), \, t \ge 0)\Big]. 
%\end{align*}

%\newpage

\subsection{Tempered stable subordinator}

Tempered $\alpha$-stable subordinators behave at short time like $\alpha$-stable subordinators and may have all moments finite,  while the latter have the first moment infinite. Therefore, we can make use of tempered stable subordinators to illustrate identities \eqref{IsoNat}--\eqref{Condition}. 
For concreteness, consider a tempered $\alpha$-stable subordinator $(\psi(t), t\ge 0)$ determined by the Laplace transform
\begin{equation} \label{TS1}
  \e e^{- u \psi(1)} = \exp\{1- (1+u)^{\alpha} \}
\end{equation}
where $\alpha\in (0,1)$. When $\alpha=1/2$, $\psi$ is also known as the inverse Gaussian subordinator. A systematic treatment of tempered $\alpha$-stable laws and processes can be found in \cite{R2}. In particular, the L\'evy measure of $\psi(1)$ is given by 
\[
\rho(dx) = \frac{1}{|\Gamma(-\alpha)|} x^{-\alpha-1} e^{-x} \, dx, \quad x>0\, ,
\]
\cite[Theorems 2.3 and 2.9(2.17)]{R2}. Therefore, the L\'evy measure $\nu$ of the process $\psi$ is given by
\begin{align}
  \nu(F) &= \int_0^{\infty}\hspace{-6pt}\int_0^{\infty} F(x \1_{[s,\infty)}) \, \rho(dx) ds \nonumber \\
  &= \frac{1}{|\Gamma(-\alpha)|} \int_0^{\infty}\hspace{-6pt}\int_0^{\infty} F(x \1_{[s,\infty)}) \, x^{-\alpha-1} e^{-x} \, dx ds\, ,  \label{LM-ts}
   \end{align}
for every measurable functional $F: \R_{+}^{[0,\infty)} \mapsto \R_{+}$.

\begin{Proposition}\label{Tempered stable} Let $ (\psi(t), t \geq 0)$ be a tempered $\alpha$-stable subordinator as above. 
\begin{itemize}
 \item[\rm{(a4)}] Given $a> 0$, let $r^{(a)}$ be the process defined by:
  \[r^{(a)}(t):= G\1_{[aU, \infty)}(t),  \quad t \geq 0\]
where $G$ has a $\mathrm{Gamma}(1-\alpha, 1)$ law and $U$ is a standard uniform random variable independent of $G$.  
Then $r^{(a)}$ satisfies \eqref{IsoNat}, that is,  
 \[
 (\psi(t) \ + \  G\1_{[aU, \infty)}(t),  \ t \geq 0)\, \=\, (\psi(t), \ t \geq 0) \ \hbox{under} \ \e\Big[\frac{\psi(a)}{\alpha a}; \ .\ \Big]
\]
  \item[\rm{(b4)}]
The L\'evy measure $\nu$ of $(\psi(t), t \ge 0)$ can be represented as  
\[
 \nu(F) = \alpha^{-1}  \e\big[G^{-1}Y e^{UY} \, F(G \1_{[UY, \infty)}) \big] 
  \]
 for every measurable functional $F: \R_{+}^{[0,\infty)} \mapsto \R_{+}$. Here $G$, $U$ are as in {\rm{(a4)}}, $Y$ is a standard exponential variable, and $G, U$ and $Y$ are independent. 
Consequently, $\nu$ is the law of the process $(G \1_{[UY, \infty)}, \,  t\ge 0)$ under the measure $ \alpha^{-1} G^{-1}Y e^{UY}  \, d\p$.

 \item[\rm{(c4)}]  The components of the decomposition \eqref{Condition}:  $\psi \= (\psi \,|\, \psi(a) =0) + {\cal L}^{(a)}$,  \  can be identified as
\[
(\psi(t), t \ge 0 \,|\, \psi(a) =0)\ \=\ (\psi(t \vee a) - \psi(a), \ t \geq 0).
\]
and 
\[
({\cal L}^{(a)}_t ,  \ t \geq 0)\ \=\ (\psi(t\wedge a), \ t \geq 0).\]

The L\'evy measures $\nu_{a}$ and $\tilde{\nu}_{a}$ of  $(\psi(t), t \ge 0 \,|\, \psi(a) =0)$  and of $({\cal L}^{(a)}_t, t \ge 0)$, respectively, are  given by
\[\nu_{a}(F) = \frac{1}{|\Gamma(-\alpha)|} \int_a^{\infty}\hspace{-6pt}\int_0^{\infty} F\big(x\1_{[s, \infty)}\big) \,  x^{-\alpha-1} e^{-x} \, dx ds \]
and
\[\tilde{\nu}_{a}(F) = \frac{1}{|\Gamma(-\alpha)|} \int_0^{a}\hspace{-6pt}\int_0^{\infty} F\big(x\1_{[s, \infty)}\big) \, x^{-\alpha-1} e^{-x} \, dx ds, \] 
  for every measurable functional $F: \R_{+}^{[0,\infty)} \mapsto \R_{+}$. 
\end{itemize}
\end{Proposition}

{\bf Proof} (a4): From \eqref{TS1} we get $\e \psi(a)= \alpha a$. Using \eqref{LM-ts}. and (\ref{r}), we get 
\begin{align*} 
  \e F(r^{(a)}_t, \, t\ge 0) &= \frac{1}{\e \psi(a)}\int F(y)\, y(a) \, \nu(dy)  \\
  &= \frac{1}{\alpha a} \int_0^{\infty}\hspace{-4pt}\int_0^{\infty} F(x \1_{[s,\infty)})\, x \1_{[s,\infty)}(a) \, \rho(dx) ds \\
  &= \frac{1}{\Gamma(1-\alpha) a} \int_0^{a}\hspace{-4pt}\int_0^{\infty} F(x \1_{[s,\infty)})\, x^{-\alpha} e^{-x} \, dx ds \\
  & =\e \big[F(G\1_{[aU, \infty)})\big]\, .
\end{align*}

\b (b4): We apply \cite[Theorem 1.2]{E} to $(\psi(t), t \geq 0)$ and $(r^{(a)}_t, t \geq 0)$ specified in (a4). Proceeding analogously to the previous examples we get for any $\sigma$-finite measure $\tilde{m}$ whose support equals $\R_{+}$ and $\int_{\R_{+}} a \, \tilde{m}(da) < \infty$
\begin{align*} 
  \nu(F) &= \int_0^{\infty} \e\left[\frac{F(r^{(a)})}{\int_0^{\infty} r^{(a)}_{s} \, \tilde{m}(ds)}\right] \e[\psi(a)]\, \tilde{m}(da) \\
  &= {1\over \alpha} \int_0^{\infty} \e\left[\frac{F(G\1_{[U_a, \infty)})}{G\tilde{m}([aU, \infty))}\right] a \, \tilde{m}(da).
\end{align*}
When $\tilde{m}$ is the law of a standard exponential  random variable we obtain
\begin{align*} 
   \nu(F) &= {1\over \alpha} \int_0^{\infty} \e\left[ e^{aU}G^{-1} F(G\1_{[aU, \infty)}) \right]  \, a  e^{-a}da \\ &= \alpha^{-1} \e\left[G^{-1}Y e^{UY} F(G \1_{[UY, \infty)}) \right]\, .
\end{align*}

(c4): We will omit this proof as it is similar to the proof of (c1) in the Poisson case.     
$\Box$

\subsection{Connection with infinitely divisible random measures}

Let ${\cal M}(S)$ denote the space of finite measures on a Borel space $(S,{\cal S})$. ${\cal M}(S)$ is a Borel space under the topology of weak convergence of finite measures. A  measurable map $\xi: \Omega \mapsto \mathcal{M}(S)$ is called a random measure on $S$. Any random measure $\xi$ can also be viewed as a stochastic process indexed by $\mathcal{S}$ and having paths in $\mathcal{M}(S) \subset \R^{\mathcal{S}}$, \ $\xi=\{\xi(A), \ A\in {\mathcal{S}} \}$. 
A random measure is called infinitely divisible if the corresponding stochastic process is infinitely divisible. 

\subsubsection{Cluster representation}
 The key result on infinitely divisible random measures is  the cluster representation. It says that any infinitely divisible random measure $\xi$ on $(S,{\cal S})$ is of the form
\begin{equation} \label{K1}
  \xi = m + \int_{{\cal M}(S)} \mu \ \Lambda(d\mu) \quad a.s.
\end{equation}
where $\Lambda$ is a Poisson random measure on ${\cal M}(S)$ with intensity $\lambda$ satisfying 
\begin{equation} \label{K2}
  \int_{{\cal M}(S)} (\mu(A) \wedge 1) \, \lambda(d\mu) < \infty, \quad A \in \mathcal{S}
\end{equation}
and $m \in {\cal M}(S)$ is non-random, see \cite[Theorem 3.20]{K}. Notice that this result follows from (\ref{NonnegativeCase}) of Section \ref{prelim} when $E = \mathcal{S}$ and $\psi= \xi$. We will sketch a proof to this claim. Indeed, since (\ref{NonnegativeCase}) in this case states that
\[ \big(\xi(A), A \in {\cal S}\big) \; \= \;  \Big(f_0(A)  \; + \,  \int_{\R_+^{\cal S}} f(A) \, N(df), \ \ A \in {\cal S}\Big), \] 
pathwise additivity of $\xi$ implies that $\nu$, the L\'evy measure of $\xi$, is concentrated on finite additive functions $f: \mathcal{S} \mapsto \R_{+}$. Since the $\sigma$-algebra $\mathcal{S}$ is countably generated and $\xi$ is pathwise $\sigma$-additive, $\nu$ is a $\sigma$-finite measure concentrated on $\mathcal{M}(S)$ with $\nu(\{0\})=0$. It follows that $f_0 \in {\cal M}(S)$.
Hence   
\[ \big(\xi(A), A \in {\cal S}\big) \; \= \;  \Big(m(A)  \; + \,  \int_{\mathcal{M}(S)} \mu(A) \, N(d\mu), \ \ A \in {\cal S}\Big). \] 
This equality can be strengthen to the almost sure equality by the usual argument. Hence \eqref{K1}-\eqref{K2} hold with $\lambda=\nu$, $\Lambda=N$, and $m=f_0$.

\subsubsection{A characterization of infinitely divisible random measures}
\label{CharacIDrm} {\rm One can make use of   (\ref{IsoNat})   for  nonnegative processes indexed by $\cal S$   to obtain the following characterization of infinitely divisible random measures on $S$.

{\sl A  random measure $\xi$ on $S$ is infinitely divisible if and only if for every $A$ in $\cal S$ such that $0 < \e[\xi(A)] < \infty$, there exists a random measure $r^{(A)}$ on $S$, independent of $\xi$ such that:
\begin{equation}\label{VeryGeneral}
\xi \ + \ r^{(A)} \  \=  \xi \ \hbox{under} \; \e[ \frac{\xi(A)}{\e[\xi(A)]},   \ . \ ]
\end{equation}

} 
The characterization (\ref{VeryGeneral}) can be connected to another characterization given in \cite[Theorem 6.17]{K}.  
  Namely, assume that $\xi$  has a $\sigma$-finite intensity $n$, then
$\xi$ is infinitely divisible if and only if for every $a$ in $S$ there exists a random measure $R^{(a)}$ on $S$, independent of $\xi$ such that
\begin{equation}\label{Olav}
\xi \, + \,   R^{(a)}  \ \=  \xi^{a},
\end{equation}
where $\xi^{a}$ is the Palm measure of $\xi$ at point $a$. 

By definition, the Palm measures $\{\xi^{a}, a \in S\}$ of $\xi$ satisfy for every $A$  in $\cal S$ and every measurable subset $L$ of ${\cal M}(S)$
\[\e[ \xi(A); \xi \in L] = \int_A n(da) \p[ \xi^a \in L],\]
which leads to the following relation for $A$ such that $0< n(A) < \infty$
\[
\p[\ \xi + r^{(A)} \in L] = \frac{1}{n(A)} \int_A n(da) \p[ \ \xi + R^{(a)} \in L]. 
\]
By computing the Laplace transforms one finally has: 
\begin{equation}
 r^{(A)}  \= \frac{1}{n(A)} \int_A n(da)  R^{(a)} . 
\end{equation}

In the special case of a point $a$ of $S$  such that $\p[\xi(\{a\}) > 0] > 0$  (e.g.  $S$ is discrete),  one obtains:  $R^{(a)} \=  r^{(\{a\})}$.

\subsubsection{A decomposition formula}

Given an infinitely divisible random measure $\xi$ on $S$, one can take advantage of (\ref{E1}) to obtain for every $A$ such that $0 < \e[\xi(A)] < \infty$, the existence of an infinitely divisible random measure ${\cal L}^{(A)}$ on $S$ such that:
\begin{equation}\label{measure}
 \xi \= (\xi \ | \ \xi(A) = 0) \; + \;   {\cal L}^{(A)},
 \end{equation}
with the two measures on the right hand side independent. }

\subsubsection{Some remarks}

In this section we take $S = \R_+^{E}$.
Let $\chi$ be a finite infinitely divisible random measure on $\R_+^{E}$ with no drift and L\'evy measure $\lambda$. 
 Assume now that for every $a$ in $E$:
\[  \int_{\R_+^{E}} f(a)  \int_{{\cal M}(\R_+^{E})} \mu(df)  \ \lambda(d\mu) < \infty .\]
 Consider then the nonnegative process $\psi$ on $E$ defined by: $\psi(x) = \int_{\R_+^{E}} f(x) \chi(df)$.  The process $\psi$ is infinitely divisible and nonnegative. The following proposition gives its L\'evy measure.

\begin{Proposition}\label{RM} The infinitely divisible nonnegative process $(\int_{\R_+^{E}} f(x) \chi(df), x \in E)$  admits for L\'evy measure $\nu$ given by:
\begin{eqnarray*} 
\nu &=&  \int_{{\cal M}(\R_+^{E})} \mu  \ \lambda(d\mu).
\end{eqnarray*}
\end{Proposition}
{\bf Proof}  From (\ref{NonnegativeCase}), we know that there exists a Poisson point process $\tilde{N}$  on $\R_+^{E}$ with intensity the L\'evy measure of $\psi$  satisfying: 
$(\psi(x), x \in E) = (\int_{\R_+^{E}} f(x) \tilde{N}(df), x \in E)$. 
Besides,    $\chi$ admits the following expression: 
$\chi =  \int_{{\cal M}(\R_+^{E})} \mu  \ N(d\mu)$, 
with  $N$  Poisson point process on ${\cal M}(\R_+^{E})$ with intensity $\lambda$.
One  obtains: 
\[(\psi(x), x \in E) = \Big( \int_{\R_+^{E}} f(x)  \int_{{\cal M}(\R_+^{E})} \mu(df)  \ N(d\mu), \quad x \in E\Big)\, . \]
Using then Campbell formula for every measurable subset $A$ of $\R_+^{E}$, one computes the intensity of the Poisson point process $ \int_{{\cal M}(\R_+^{E})} \mu(df)  \ N(d\mu)$
\[ \e[ \int_{\R_+^{E}} 1_A(f)  \int_{{\cal M}(\R_+^{E})} \mu(df)  \ N(d\mu)] =  \int_{\R_+^{E}} 1_A(f)  \int_{{\cal M}(\R_+^{E})} \mu(df)  \ \lambda(d\mu)= \nu(A)\, . \]
$\Box$

Proposition \ref{RM} allows to write every L\'evy measure $\nu$ on $\R^E_+$ in terms of a L\'evy measure on ${\cal M}(\R_+^{E})$. Indeed given a L\'evy measure $\nu$ on $\R^E_+$, denote by  $(\psi(x), x \in E)$ the corresponding infinitely divisible nonnegative process without drift.  From (\ref{NonnegativeCase}), we know that $\psi$ admits the representation $(\int_{\R_+^{E}} f(x) \chi(df), x \in E)$ with $\chi$ Poisson random measure on ${\cal M}(\R_+^{E})$. The random measure $\chi$ is hence infinitely divisible. Proposition \ref{RM} gives us:  
\begin{equation}\label{notunique}
\nu = \int_{{\cal M}(\R_+^{E})} \mu  \ \lambda(d\mu)
\end{equation}
 where $\lambda$ is the L\'evy measure of $\chi$. Proposition \ref{RM} allows to see that given $\nu$,    the L\'evy measure $\lambda$ satisfying (\ref{notunique})  is not unique.

\begin{Proposition}\label{poisson}
The intensity $\nu$ of a  Poisson random measure on $\R_+^E$  with L\'evy measure $\lambda$  satisfies:
\[ \nu = \int_{{\cal M}(\R^E_+)} \mu \ \lambda(d\mu).\]
\end{Proposition}

\subsection{Infinitely divisible permanental processes}

\b A permanental process $(\psi(x) , x \in E)$ with index $\beta > 0$ and kernel $k =(k(x,y), (x,y) \in E\times E)$
is a nonnegative process with finite dimensional Laplace transforms satisfying, for every $\alpha_1,.., \alpha_n \geq 0$ and  every $x_1$, $x_2$,..,$x_n$ in $E$:
\begin{equation}\label{perm}
\e[ \exp\{- {1\over 2}\sum_{i = 1}^n \alpha_i \psi(x_i) \}] = \det( I + \alpha K)^{-\beta}
\end{equation}
where $\alpha$ is the diagonal matrix with diagonal entries $(\alpha_i)_{1\leq i \leq n}$, $I$ is the $n\times n$-identity matrix and $K$ is the matrix $(k(x_i,x_j))_{1 \leq i,j \leq n}$. 
\smallskip

\b Note that the kernel of a permanental process is not unique. 

\smallskip

\noindent In case $\beta = 1/2$ and $k$ can be chosen symmetric positive semi-definite, $(\psi(x), x \in E)$  equals in law $(\eta^2_x, x \in E)$ where $(\eta_x , x \in E)$ is a centered Gaussian process with covariance $k$.  The permanental processes hence represent an extension of the definition of squared Gaussian processes.

A necessary and sufficient condition on $(\beta, k)$ for the existence of a permanental process $(\psi(x), x \in E)$ satisfying (\ref{perm}),  has been established by Vere-Jones \cite{VJ}.  
 Since we are interested by the subclass of infinitely divisible permanental processes, we will only remind a necessary and sufficient condition for a permanental process to be infinitely divisible.  Remark that if $(\psi(x), x \in E)$ is infinitely divisible then for every measurable nonnegative $d$, $(d(x)\psi(x), x \in E)$ is also infinitely divisible.  Up to the product by a deterministic function, 
 $(\psi(x), x \in E)$ is infinitely divisible if and only if it admits for kernel the $0$-potential densities (the Green function)  of a transient Markov process on $E$ (see \cite{EK} and \cite{EM}). 

\b Consider an infinitely divisible permanental process $(\psi(x), x \in E)$ admitting for kernel the Green function $(g(x,y), (x,y) \in E\times E)$ of a transient Markov process $(X_t, t \geq 0)$ on $E$.  
For simplicity assume that $\psi$ has index $\beta = 1$.  For $a \in E$ such that $g(a,a) > 0$,  denote by $(L^{(a)}_{\infty}(x), x \in E)$ the total accumulated local times process of $X$ conditioned to start at $a$ and killed at its last visit to $a$.  In \cite{E}, 
 (\ref{Condition}) has been explicitly written for $\psi$:
\[ \psi \=  (\psi | \psi(a) = 0) \ + \ {\cal L}^{(a)} \]
with  ${\cal L}^{(a)}$ independent process of $(\psi | \psi(a) = 0)$, such that
${\cal L}^{(a)} \=  (2 L^{(a)}_{\infty}(x), x \in E)$.  Moreover $(\psi | \psi(a) = 0)$  is a permanental process with index $1$ and with kernel the Green function of $X$ killed at its first visit to $a$. 

One can also explicitly write (\ref{IsoNat}) for $\psi$ with $(r^{(a)}(x), x \in E) \= (2 L^{(a)}_{\infty}(x), x \in E)$.  Hence the case of infinitely divisible permanental processes is a special case since $r^{(a)}$ is  infinitely divisible and $r^{(a)} \= {\cal L}^{(a)}$. 

The easiest way to obtain the L\'evy measure $\nu$ of $\psi$ is to use (\ref{levyMeasure}) with $m$ $\sigma$-measure with  support equal to $E$ such that: $\int_E g(x,x) m(dx) < \infty$,  to obtain
\[\nu(F) =  \int_E \e[\frac{F(2 L_{\infty}^{(a)})}{\int_E L^{(a)}_{\infty}(x) m(dx)}] g(a,a) m(da),\]
for any measurable functional $F$ on $\R^{E}_+$.

\b  If moreover,  the $0$-potential densities  $(g(x,y), (x,y) \in E\times E)$ were taken with respect to $m$ then,  for every $a$,
$ \int_E L^{(a)}_{\infty}(x) m(dx)$ represents the time of the last visit to $a$ by $X$ starting from $a$. 

\bigskip

\section{Transfer of continuity properties}\label{Correspondency}
 
Using \eqref{E1}, a nonnegative infinitely divisible process $\psi = (\psi(x), x \in E)$ with L\'evy measure $\nu$ and no drift, is hence connected to a family of nonnegative infinitely divisible processes $\{ {\cal L}^{(a)}, a \in E\}$.  In case when $\psi$ is an infinitely divisible squared Gaussian process, Marcus and Rosen \cite{MR} have established correspondences between path properties of $\psi$ and the ones of ${\cal L}^{(a)}, a \in E$.  
To initiate a similar study for a general $\psi$, we 
assume that $(E,d)$ is a separable metric space with a dense set $D = \{a_k , k \in \N^*\}$. 

One immediately notes that if $\psi$ is a.s. continuous with respect to  $d$, then for every $a$ in $E$, ${\cal L}^{(a)}$ is a.s. continuous  with respect to  $d$ and the measure $\nu$ is supported by the continuous functions from $E$ into $\R_+$ i.e. $r^{(a)}$ is continuous with respect to  $d$ , for every $a$ in $E$. 

Conversely if ${\cal L}^{(a)}$ is continuous with respect to  $d$ for every $a$ in $E$, what can be said about the continuity of $\psi$ ?

As noticed in \cite{R} (Proposition 4.7) the measure $\nu$ admits the following decomposition:
\begin{equation}\label{levyDecomp}
\nu \ = \ \sum_{k = 1}^{\infty}  1_{A_k} \nu_k,
\end{equation}
where $A_1 =  \{ y \in \R^E_+ :   y(a_1) > 0\}$ and for $k > 1$, 
\[A_k = \{ y \in \R^E_+ :  y(a_i) = 0, \forall i < k  \ \hbox{and} \  y(a_k) > 0\}\]  
and $\nu_k$ is defined by
\[\nu_k(F)  = \e[ \frac{\e(\psi(a_k)}{r^{(a_k)}_{a_k}} 1_{A_k}(r^{(a_k)}) F(r^{(a_k)})] \]
for every measurable functional $F: \R_{+}^{E} \mapsto \R_{+}$.

For every $k$ the measure $\nu_k$ is a L\'evy measure.  Since the supports of this measures are disjoint they correspond to independent nonnegative  infinitely divisible processes that we denote by $ L(k)$, $k\geq 1$. As a consequence of (\ref{levyDecomp}), $\psi$ admits the following decomposition:
\begin{equation}\label{IDdecomp}
\psi \ \= \ \sum_{k = 1}^{\infty} L(k).
\end{equation}

 Note that 
\[L(1) \=  {\cal L}^{(a_1)}\]
and similarly for every $k > 1$:
\[ L(k)  \= ({\cal L}^{(a_k)} | {\cal L}^{(a_k)}_{|_{\{a_1,..,a_{k-1}\}}} = 0).\]
Consequently, for every $k \geq 1$, $L(k)$ is continuous with respect to  $d$. 

From (\ref{IDdecomp}), one obtains all kind of $0-1$ laws for $\psi$. For example:

\b - $\p[\psi \ \hbox{is continuous on}\  E] = 0$ or $1$.

\b - $\psi$ has a deterministic oscillation function $w$, such that for every $a$ in $E$:
\[ \liminf_{x \rightarrow a} \psi(x) =\psi(a) \ \hbox{and} \ \limsup_{x \rightarrow a} \psi(x) = \psi(a) + w(a).\]

\smallskip

\b Exactly as in \cite{E}, one shows the following propositions.

\begin{Proposition} If for every $a$ in $E$, ${\cal L}^{(a)}$ is a.s. continuous, then there exists a dense subset $\Delta$ of $E$ such that a.s. $\psi$ is continuous at each point of $\Delta$ and $\psi_{|_{\Delta}} $ is continuous.
\end{Proposition}

\begin{Proposition} Assume that $\psi$ is stationary. If for every $a$ in $E$, ${\cal L}^{(a)}$ is a.s. continuous, then $\psi$ is continuous. 
\end{Proposition}

\section{A limit theorem}\label{limit}
Given a nonnegative infinitely divisible without drift  process $(\psi_x, x \in E)$, the following result gives an intrinsic way to obtain $r^{(a)}$ for every $a$ in $E$. 

\begin{Theorem}\label{natural} For a nonnegative infinitely divisible process $(\psi_x , x \in E)$ with L\'evy measure $\nu$,  denote by $\psi^{(\delta)}$ an infinitely divisible process with L\'evy measure $\delta \nu$.  Then, for any $a$ in $E$ such that $\e[\psi_a] > 0$, $r^{(a)}$ is the limit in law of the processes $\psi^{(\delta)}$ under $\e\Big[\frac{\psi^{(\delta)}_a}{\e[\psi^{(\delta)}_a]}; \, \cdot \, \Big]$, as $\delta \to 0$. 
\end{Theorem}

\b {\bf Proof}  We remind (\ref{r}): 
$\p[r^{(a)} \in dy] = \frac{y(a)}{\e[\psi_a]} \, \nu(dy)$. 
Since $\e[\psi^{(\delta)}_a] = \delta \e[\psi_a]$, one obtains immediately:  $\p[r^{(a)} \in dy] = \frac{y(a)}{\e[\psi^{(\delta)}_a]}  \ \delta\nu(dy)$. Consequently $r^{(a)}$ satisfies: 
\[ \psi^{(\delta)} \ + \ r^{(a)} \ \= \ \psi^{(\delta)}  \ \hbox{under} \  \e[\frac{\psi^{(\delta)}_a}{\e[\psi^{(\delta)}_a]}; \ . \ ] .\]
As $\delta \to 0$,  $\psi^{(\delta)}$ converges to the $0$-process in law, so 
$\psi^{(\delta)}$ under $\e[\frac{\psi^{(\delta)}_a}{\e[\psi^{(\delta)}_a]}; \, \cdot  \, ]$ must converge in law to $r^{(a)}$.  $\Box$

\bigskip
\b From (\ref{IsoNat}) and (\ref{Condition}), one obtains in particular: 
\begin{equation}
{\cal L}^{(a)} \ + \ r^{(a)} \=  {\cal L}^{(a)} \, \hbox{under} \; \e[\frac{{\cal L}^{(a)}_a}{\e[{\cal L}^{(a)}_a]}; \ . ]
\end{equation}
We know from \cite{E}, that the L\'evy measure of ${\cal L}^{(a)}$ is $\nu(dy) 1_{y(a) > 0}$. 
Denote by $\l^{(a, \delta)}$ a nonnegative process with L\'evy measure $\delta \nu(dy) 1_{y(a) > 0}$.  Using Theorem \ref{natural}, one obtains that $r^{(a)}$ is  also the limit in law of ${\l}^{(a,\delta)} \  \ \ \hbox{under} \;  \e[ \frac{ {\l}^{(a,\delta)}_a}{\e[{\l}^{(a,\delta)}_a]}; \ . \ ]$.

\end{document}